\documentclass[reqno, 12pt]{amsart}

\def\C{{\mathbb C}}

\def\Q{{\mathbb Q}}
\def\Qp{{\mathbb Q}_p}
\def\Z{{\mathbb Z}}
\def\Zp{{\mathbb Z}_p}
\def\R{{\mathbb R}}
\def\A{{\mathbf A}}
\def\K{{\mathbf K}}
\def\k{{\mathbf k}}
\def\bfC{{\mathbf C}}

\def\Aut{\operatorname{Aut}}
\def\Gal{\operatorname{Gal}}

\def\Qpbar{{\overline{\Q}_p}}
\def\Fp{{{\mathbb F}_p}}
\def\F{{\mathbb F}}
\def\Qp {{{\mathbb Q}_p}}
\def\Fq{{{\mathbb F}_q}}
\def\Fr{{{\mathbb F}_r}}

\newtheorem{lemma}{Lemma}
\newtheorem{cor}{Corollary}
\newtheorem{prop}{Proposition}

\newtheorem{theorem}{Theorem}

\theoremstyle{definition}
\newtheorem{defn}{Definition}   

\newtheorem{question}{Question}

\theoremstyle{remark}
\newtheorem{rem}{Remark}        
      
\newtheorem{example}{Example}

\begin{document}
\title[The impact of the infinite primes]
{The impact of the infinite primes on the Riemann Hypothesis for 
Characteristic $p$ $L$-series}
\author{David Goss}
\thanks{This paper is dedicated to Ram with great respect and affection on
his $70$-th birthday}
\address{Department of Mathematics\\ The Ohio State University\\ 231 W.
$18^{\text{th}}$ Ave. \\ Columbus, Ohio 43210}
\email{goss@math.ohio-state.edu}
\date{May 25, 2001}

\begin{abstract}
In \cite{go2} we proposed an analog of the classical Riemann hypothesis
for characteristic $p$ valued $L$-series based on the
work of Wan, Diaz-Vargas, Thakur, Poonen, and Sheats for the zeta
function $\zeta_{\Fr[\theta]}(s)$.
During the writing of \cite{go2}, we made
two assumptions that have subsequently proved to be incorrect. The first 
assumption is that we can ignore the trivial zeroes of characteristic $p$
$L$-series in formulating our conjectures. Instead, we show here how 
the trivial zeroes influence nearby zeroes and so lead to counter-examples
of the original Riemann hypothesis analog.
We then sketch an approach to handling such ``near-trivial'' zeroes via
Hensel's and Krasner's Lemmas (whereas classically one uses
Gamma-factors). Moreover, we show that 
$\zeta_{\Fr[\theta]}(s)$ is not representative of
general $L$-series as, surprisingly, all its zeroes
are near-trivial, much as the Artin-Weil zeta-function of $\mathbb{P}^1/\Fr$
is not representative of general complex $L$-functions of curves over finite
fields. Consequently, the ``critical zeroes'' 
(= all zeroes not effected by the trivial zeroes) of characteristic
$p$ $L$-series now appear to be quite mysterious.
The second assumption made while writing
\cite{go2} is that certain Taylor expansions of classical
$L$-series of number fields
would exhibit complicated behavior with respect to their
zeroes. We present a simple argument that this is not so, and, at the same time,
give a characterization of functional equations.
\end{abstract}

\maketitle

\section{Introduction}\label{intro}
In the paper \cite{go2} (whose notations etc., we generally follow here)
we defined a possible Riemann hypothesis for the zeroes of
characteristic $p$ $L$-functions. This work is based on the results of Wan 
\cite{w1}, Diaz-Vargas and
Thakur \cite{dv1}, Poonen, and Sheats \cite{sh1} for the zeta
function $\zeta_{\Fr[\theta]}(s)$ ($s=(x,y)\in S_\infty$; see
Section \ref{trivial} for the definitions). In particular
the zeroes of $\zeta_{\Fr[\theta]}(s)$ were found to be both
simple and lie on the line $\Fr((1/T))$ just as the completed Riemann zeta
function $\pi^{-s/2}\Gamma(s/2)\zeta(s)$,
$s=1/2+it$, is conjectured to have only real zeroes in $t$.
This finite characteristic Riemann hypothesis has two components
to it. The first, Conjecture 4 of \cite{go2}, 
is a basic finiteness statement
about the extension fields generated by $L$-series zeroes (obviously
no such conjecture need be made classically). The second,
Conjecture 5 of \cite{go2}, focuses on
counting the number of zeroes of a given absolute value; based on
the above examples, we postulated that, outside of finitely many exceptions,
the absolute value determined the zero. (Both conjectures are recalled 
in Section \ref{review} below; therefore we will drop the reference
to \cite{go2} from now on.)
We pointed out similarities between Conjecture 5 and the classical
Riemann hypothesis for number fields.

Now in 
classical theory we work, of course, with the completed
$L$-function including the Gamma-factors arising from the infinite
primes as above for $\zeta(s)$. These Gamma-factors have poles at
some of the negative integers (e.g., $\Gamma(s/2)$ has 
poles at the non-positive
even integers) and so, by analytic continuation, they force the $L$-series
to have associated ``trivial zeroes.'' The functional equation then
assures us that all other zeroes must lie in the critical strip.
On the other hand, in the characteristic $p$ case presented
in \cite{go2} we ignored 
input from the infinite primes and their associated trivial
zeroes. Indeed, our conjectures allow for finitely
many exceptional cases for each $y$,
and, for $y$ a negative integer, there are only
finitely many trivial zeroes; thus ignoring them appeared innocuous.
However, since the publication of \cite{go2}, we have realized
that the infinite primes must also be taken into account in the
characteristic $p$ Riemann hypothesis. It is our goal here to explain
how the trivial zeroes lead to counter-examples for Conjecture 5
(though Conjecture 4 still appears to
be valid). More precisely, the topology of our space $S_\infty$, which
is the
natural domain of definition for the characteristic $p$ $L$-series, permits
us to inductively construct counter-examples using zeroes sufficiently
close to trivial zeroes (we call such zeroes ``near-trivial zeroes'').
We also suggest
some possible replacements for Conjecture 4 by sketching a procedure
based on Hensel's Lemma to isolate the near-trivial zeroes.
The original conjectures may then
work for the remaining zeroes (which, following
classical precedent, we call the ``critical zeroes'').

However, there is now a great surprise. Upon searching for
the critical zeroes of $\zeta_{\Fr[\theta]}(s)$, $s\in S_\infty$, we find
that they do not exist! (More precisely, we prove this fact for
$\Fp[\theta]$; the general $\Fr[\theta]$ case seems very
likely to follow from Sheats' techniques --- see also
Corollary \ref{evenval} of
Section \ref{basic}.) That is, the case of $\zeta_{\Fr[\theta]}(s)$
may not be representative of general characteristic
$p$ $L$-series, much as the Artin-Weil zeta-function of
$\mathbb{P}^1$ over a finite field is not typical of the $L$-series of
general curves. Indeed, the Artin-Weil zeta-function of the projective
line also has no critical zeroes. We will exhibit here
a few examples of critical
zeroes in the characteristic $p$ theory and we hope to have more examples
in the near future. However, as of this
writing, they are a complete mystery. For instance where they lie, or even
if there are infinitely many of them (for a given $L$-series and
interpolation place) is not known. It is certainly possible that there
may be more surprises ahead.

Let $K$ be a complete, algebraically closed, non-Archimedean field 
and let $f(x)=\sum a_nx^n$ be an entire power series with coefficients in
$K$. Let $L$ be a complete subfield of $K$ which contains the zeroes of
$f(x)$. It is then a standard fact that there exists 
$\alpha\in K^\ast$ such that the coefficients of $\alpha f(x)$ lie
in $L$; i.e., the coefficients of $f(x)$ lie on an $L$-line through the
origin in $K$. As is well-known, 
in complex analysis the relationship between the zeroes
and the Taylor coefficients of an entire function is more complicated;
e.g., $e^{2\pi i x}-1$. 
(However, as Keith Conrad pointed out, a similar statement is true in
complex analysis if we allow ourselves to multiply by an 
invertible entire function; in non-Archimedean analysis such entire functions
are constant.)
During the writing of \cite{go2}, we had assumed that
Taylor expansions of classical $L$-series 
(of number fields) would exhibit similar complicated behavior. 
This assumption has also proved to be incorrect. Indeed, in our last section
we give a simple argument which shows that the associated Taylor coefficients
in $t$ (where $s=1/2+it$) lie on a line through the origin in $\mathbb C$.
This argument
also characterizes the functional equation of the $L$-series
via the description of complex conjugation given in Equation \ref{funceq}.

We now describe in more detail the contents of this paper. In Section
\ref{trivial} we review the definitions and interpolations
of characteristic $p$ $L$-series. We recall how  the analogue
of Artin $L$-series has trivial zeroes (essentially from the classical
theory of Artin and Weil) at the negative integers. 
We also present a reasonable approach to trivial zeroes for the general
$L$-series of a Drinfeld module; it is our belief that the work of
Boeckle and Pink \cite{bp1}, \cite{b1} will ultimately flesh out
this construction.

In Section \ref{krasner} we recall Krasner's Lemma in order to put it into
a form more useful for calculations in characteristic $p$. In Section
\ref{review} we review the statements of Conjectures 4 and 5.

In Section \ref{counter} we construct counter-examples to Conjecture
5. Along the way we establish that the family of Newton
polygons associated to an $L$-series possesses certain invariance
properties, see Lemma \ref{newtonclose}. More precisely, the first
$n$-segments of the Newton polygon will, generically, be an invariant
of a natural group of translations. This leads us to believe that the family of
Newton polygons itself will serve to distinguish between $L$-series
(see Question \ref{langlands} in this section). We also sketch
the beginnings of an approach to remove these counter-examples by
isolating the offending near-trivial zeroes via Hensel's Lemma. We give some
examples of how this will work (as well as find a few critical zeroes)
but total success here will need many more results on the structure
of the zeroes.

In Section \ref{basic} we use the techniques of Diaz-Vargas and
Sheats to study the effects of the trivial zeroes for $\zeta_{\Fr[\theta]}(s)$,
$s\in S_\infty$, and we establish in Proposition \ref{jandn}
that all zeroes for $\zeta_{\Fp[\theta]}(s)$ are near-trivial. It is our
belief that ultimately there should be some massive generalization of
the techniques of Diaz-Vargas and Sheats
to general characteristic $p$ $L$-series enabling us, at least, to 
specify exactly where the near-trivial zeroes lie. Indeed, Boeckle
has established the first of such results \cite{b1} by showing the
logarithmic growth of the degrees of special polynomials in general. For
$\Fr[T]$, this logarithmic growth is a first corollary of the techniques
of \cite{dv1} and \cite{sh1}. As one sees by simple examples (such
as in the discussion directly after Corollary \ref{log}), this logarithmic
growth is analogous to the formula giving the growth of the number of
zeroes of classical $L$-series in the critical strip. 

Finally, let $L(\chi,s)$, $s=1/2+it\in \C$, be the
$L$-series associated to a number field
and abelian character $\chi$. Let $\Lambda(\chi,s)$ be the
completed $L$-series which includes the Gamma-factors. In our last
section, Section \ref{taylor}, we present an elementary
argument that the Taylor coefficients of
$\Lambda(\chi,s)$ about $t=0$ are, up to multiplication be a non-zero constant, real
numbers. This is contrary to what was expected during the writing of
\cite{go2}.

It is my pleasure to thank Gebhard Boeckle, Keith Conrad, Mike Rosen,
and Dinesh Thakur for their comments on an earlier
version of this manuscript.
\section{Trivial zeroes} \label{trivial}
As mentioned in the introduction,
we will follow the notation of \cite{go2}. We begin by briefly reviewing the
definition of $L$-functions of Galois characters, Drinfeld modules,
etc. Let $\mathcal X$ be a smooth, projective, geometrically
connected curve over the finite field $\Fr$ where $r=p^m$ and $p$ is 
prime. Let $\infty\in \mathcal X$ be  fixed closed point of degree
$d_\infty$ over $\Fr$.
Let $\k$ be the function field of $\mathcal X$ and let $\A$ be the
subring of those functions which are regular outside of $\infty$. In the
theory, one views $\k$ as the analogue of $\Q$ and $\A$ as the analogue of
$\Z$; indeed, $\A$ is a Dedekind domain with finite class group and
group of units $\Fr^\ast$. 

Let $w$ be an arbitrary place of $\k$. We let $\vert x\vert_w$ be the
normalized absolute value (= multiplicative valuation)
at $w$ with associated additive valuation
$\nu_w(x)$; by definition $\nu_w(t)=1$ if $t\in \k$ has a simple
zero at $w$. We let $\k_w$ be the completion of $\k$ with respect to
$w$ and we denote its finite field of constants by $\F_w$. We let
$d_w:=[\F_w\colon \Fr]$ be the degree of $w$ over $\Fr$.
We let $\bfC_w$ be the completion of a fixed algebraic closure
of $\k_w$ equipped with the canonical extension of $\vert x\vert_w$.
In particular we put $\K:=\k_\infty$.
Thus $\K, {\mathbf C}_\infty$ are the analogues of
$\R, \C$, respectively, whereas $\k_v$, for a finite $v$, is the analogue
of $\Qp$ etc.

We now let $A,k,K$ etc., be another copy of these rings. There is
an obvious isomorphism $\theta$ from the ``bold'' to the ``non-bold'' which
makes $k,K, C_\infty$ into $\A$-fields. As in \cite{go2}, we view the non-bold
fields as being the ``scalars'' over which we can define  Drinfeld modules
etc., equipped with an action by the ``operators'' in $\A$. 
The basic example is the Carlitz module $C$ defined
for $\A=\Fr[T]$. Here we put $\theta=\theta(T)\in k$ and 
$$C_T:=\theta \tau^0+\tau,$$
where $\tau\colon C_\infty\to C_\infty$ is the $r$-th power mapping of the
field $C_\infty$; thus
the Carlitz module is obviously defined over $k=\Fr(\theta)$. In fact, it is
easy to see that $C$ 
gives rise to an obvious family of Drinfeld modules over ${\rm Spec}(A)$.
(It will always be clear to the reader when ``$C$'' is being used to denote
a field or the Carlitz module.)

More generally, let $L\subset C_\infty$ 
be a finite extension of $k$ and let $\omega$ be a place
of $L$. We say $\omega$ is a ``finite place'' (or finite prime)
if it lies over a prime of
$A$; otherwise it is an ``infinite place'' (or infinite prime).
This notion carries over
to $\k$ etc., in the obvious fashion.

Let $\psi$ be
a Drinfeld $\A$-module of rank $d$ defined over $L$. 
Let $\mathfrak P$ be a finite prime of $L$ lying over the prime
$\mathfrak p$ of $A$ with finite residue class fields $\F_\mathfrak P$ and
$\F_\mathfrak p$ respectively.
As usual, we define the norm
$n\mathfrak P$ of $\mathfrak P$ to
be ${\mathfrak p}^f$ where $f$ is the residue field degree.
All but finitely many such primes $\mathfrak P$
are good for $\psi$ in that one can reduce $\psi$ modulo $\mathfrak P$ to 
obtain a Drinfeld module $\psi^\mathfrak P$ over
$\mathbb{F}_\mathfrak P$ with the same rank
as $\psi$. Associated to $\mathbb{F}_\mathfrak P$ there is the Frobenius
endomorphism ${\rm Fr}_\mathfrak P$ of $\psi^\mathfrak P$, and one sets
$$P_\mathfrak P(u):=\det(1-u{\rm Fr}_{\mathfrak P}
\mid T_v(\psi^{\mathfrak P}))\,;$$
here $v\subset \A$ is a non-trivial prime such that
$\theta(v)\neq \mathfrak p$ and
$T_v$ is the $v$-adic {\it Tate module} of $\psi^\mathfrak P$. It is easy to
see that the
$\A_v$-module
$T_v(\psi)$ is free of
rank $d$. In complete
accordance with classical theory, $P_{\mathfrak P}(u)$ has $\A$-coefficients
which are independent of the choice of $v$. Moreover, as Drinfeld
has shown, its roots inside $\bfC_\infty$
satisfy the local Riemann hypothesis in terms of their
absolute values, see, e.g., \S 4.12 of \cite{go1}.

Now assume that we have chosen a notion of ``sign'' on $\K$; that is,
a homomorphism ${\rm sgn}\colon \K^\ast \to \F_\infty^\ast$ which is
the identity on $\F_\infty^\ast\subset \K^\ast$. Elements $x$
with ${\rm sgn}(x)=1$ are called ``positive.'' Let $\pi\in \K^\ast$
be a fixed positive uniformizer and
let $a\in\K$ be a positive element with a pole of order $e$ at $\infty$. We
set 
$$\langle a \rangle=\langle a \rangle_\pi:=\pi^e a\,.$$
It is clear that $\langle a \rangle\in U_1$, where
$U_1$ is the multiplicative group of $1$-units in $\K^\ast$, and that
$\langle ab\rangle=\langle a\rangle \langle b \rangle$ for
positive $a$ and $b$. Note that the binomial theorem makes
$U_1$ into a $\Zp$-module.

Let $a$ be a positive element of $\A$. We set
$\langle (a)\rangle:=\langle a \rangle$ thereby giving a homomorphism from
the group of principal and positively generated $\A$-fractional ideals
to $U_1$. Let $\hat{U}_1 \supset U_1$ be the group of all $1$-units in $\bfC_\infty$.
As $\hat{U}_1$ may be seen to be a $\Qp$-vector space ($p$-th roots
may be uniquely taken),
a simple argument implies that $\langle ?\rangle$
extends uniquely to a homomorphism from the group of all fractional ideals
to $\hat{U}_1$.
Let $\K_{\mathbf V}\subset \bfC_\infty$ be the subfield obtained by
adjoining $\langle I\rangle$ to $\K$ where $I$ ranges over all
fractional ideals of $\A$. As $U_1$ is a $\Zp$-module
and the class number of $\A$ is finite, we conclude that $\K_{\mathbf V}$
is a finite, totally inseparable, extension of $\K$.

Let $S_\infty:=\bfC^\ast_\infty\times
\Zp$ and let $s=(x,y)\in S_\infty$. For a non-zero ideal $I\subseteq \A$
we put
$$I^s:=x^{\deg I}\langle I\rangle^y\,.$$
Let $\pi_\ast\in \bfC_\infty$ be a fixed $d_\infty$-th root of $\pi$
and let $j$ be an integer. It is then easy to see that
$$(i)^{s_j}=i^j$$
for positive $i\in \A$ and $s_j:=(\pi_\ast^{-j},j)\in S_\infty$.
We frequently write ``$j$'' for ``$s_j$.''

By abuse of language, we also write $n\mathfrak P$ for
$\theta^{-1}(n\mathfrak P)$ whenever no confusion will
arise. Thus, finally, the
$L$-series $L(\psi,s)$, $s\in S_\infty$, of $\psi$ over the field
$L$ is defined by
$$L_\infty(\psi,s)=L(\psi,s):=
\prod_{{\mathfrak P}~\rm good}P_{\mathfrak P}(\theta^{-1}(n{\mathfrak P })^{-s})^{-1}=
\prod_{{\mathfrak P}~\rm good}P_{\mathfrak P}(n{\mathfrak P }^{-s})^{-1}\,.$$
The local Riemann hypothesis assures us that these $L$-series converge
on the ``half-plane'' of $S_\infty$ defined by
$\{(x,y)\mid \vert x \vert_\infty
> t\}$, for some positive real $t$, to a ${\mathbf C}_\infty$-valued
function. In a similar fashion one can construct
$L$-series of pure $T$-modules etc.

Ultimately, as with classical theory,
Euler factors at the finitely many bad primes 
ought to be added into
the definition, see Remark \ref{whybadfactors}.
Such factors should be given by recent work of F.\ Gardeyn \cite{ga1}.
In any case, unlike classical theory, there are also many
examples where all finite primes are good. 

Now let $L^{\rm sep} \subset C_\infty$ be the separable closure of $L$. Let
$\Qpbar$ be a fixed algebraic closure of $\Qp$ and let
$V$ be a finite dimensional $\Qpbar$-vector space. Let
$\rho\colon \Gal (L^{\rm sep}/L)\to \Aut_{\Qpbar} (V)$ be a 
representation of Galois type (i.e., which factors through the Galois group
of a finite Galois extension of $L$). As explained in 
\S 8 of \cite{go1}, the classical
definition of Artin $L$-series is easily modified to define a
${\mathbf C}_\infty$-valued $L$-series
$L(\rho,s)$, $s\in S_\infty$, whose Euler product converges on the half-plane
$\{(x,y)\mid \vert x\vert_\infty>1\}$.

We will refer to $L(\psi,s)$, for $\psi$ a Drinfeld module, as 
an $L$-series of ``Drinfeld type,'' and $L(\rho,s)$, for
a Galois representation $\rho$, as an $L$-series of ``Galois
type.''
Both types of $L$-series will be referred to as
``$L$-series of arithmetic type'' from now on. 

Let $L(s)$, 
$s=(x,y)\in S_\infty$, be an $L$-series of
arithmetic type. For fixed $y\in \Zp$, $L(x,y)$
is an power series in $x^{-1}$ with coefficients in a
finite extension $\K_L(y)$ of $\K$. (To avoid any possibility of
confusion we will always use a subscripted ``$L$'' to refer to 
a particular $L$-series and not a field $L$.)
For instance, if $L(s)=L(\psi,s)$ is an
$L$-series of Drinfeld type, then $\K_L(y)$ is a subfield of
$\K_{\bf V}$. In fact, let $u$ be the order of the strict class group of $\A$
(which is the quotient of the group of all
$\A$-fractional ideals modulo the subgroup of principal and positively generated
ideals). Then if $y\in u\Zp$, one sees easily that $\K_L(y)=\K$. If 
$L(s)$ is an $L$-series of Galois type, then $\K_L(y)/\K$ may contain a finite extension
of constant fields.

Let $\pi_\ast\in \bfC_\infty$ be a fixed $d_\infty$-th root of our
fixed parameter $\pi$,
and let $j$ be a non-negative integer.
As is standard, we put
$$z_L(x,-j):=L(x\pi^j_\ast,-j)\,.$$
The finiteness of the class number of $\A$ implies that $z_L(x,-j)$ is
a power series whose coefficients lie in a finite extension of
$\k$ and are integral over $\A$.
For $L(\psi,s)$ of Drinfeld type,
recent work \cite{bp1}, \cite{b1} expresses these power series in terms of the
cohomology of certain ``crystals''  and thus establishes that they are
actually polynomials in $x^{-1}$. A similar statement in the case
of an $L$-series of Galois type had previously been shown using elementary
estimates and the classical theory of Weil.

\begin{defn}\label{specialpolys}
The polynomials $z_L(x,-j)$ are the {\it special polynomials} of $L(s)$.
\end{defn}

The
cohomological description  of the special polynomials is critical in 
the analytic continuation of an $L$-series  
$L(\psi,s)$ of Drinfeld type to an essentially algebraic entire function on all of
$S_\infty$ (as defined in \S 8.5 of \cite{go1}). Indeed, as in \cite{b1}, one
may use the cohomological description to
give a {\it logarithmic} bound on the growth of the degrees (in
$x^{-1}$) of $z_L(x,-j)$ as a function of $j$ from which
the analytic continuation is readily deduced. For the
Carlitz module $C$, whose $L$-series is easily seen to be
$\zeta_{\Fr[\theta]}(s-1)$, this bound was originally
shown by H.\ Lee \cite{le1} using elementary methods; see Th.\ VIII of \cite{th1} for a statement
of related results.
This bound also follows from the work of Diaz-Vargas and Sheats
as in Section \ref{basic}. In fact, the logarithmic bound on the degrees of
the special polynomials for any $L$-series of Galois type
may be established by elementary, non-cohomological means.

Let $L(\rho,s)$ be an $L$-series of
Galois type. In Sections 8.12 and 8.17 of \cite{go1}
there is a ``double congruence'' relating
$z_L(x,-j)$ to the incomplete characteristic $0$-valued $L$-function
$\hat{L}(\rho\otimes \omega^{-j},t)$ of $\rho$ 
twisted by powers of the Teichm\"uller character; this incomplete $L$-series is
defined by the usual Euler product but taken only over the
finite primes. 
By Weil's Theorem (= the Artin Conjecture in this context), this
incomplete $L$-series is a polynomial in $t=r^{-s}$ which is divisible by
the finite Euler product taken only over the infinite primes.
By using such double congruences infinitely often, one
deduces that $z_L(x,-j)$, and thus $L(\rho,(x,-j))$, have a number of canonical
zeroes.

\begin{defn} \label{trivialzero}
The zeroes just described are the called the trivial zeroes 
of $L(\rho,s)$, and $z_L(x-j)$, at $y=-j$.
\end{defn}

\noindent
We will use the expression ``trivial zeroes'' to refer to the union
of the trivial zeroes at $-j$ for all $j$. It is very easy to see that
the trivial zeroes belong to a finite extension of $\K$.

Note that as $L(\rho, (x,-j))$ is a polynomial
in $x^{-1}$, there are obviously only finitely many trivial zeroes for
each $j$. As $\rho$ is of Galois type, the trivial zeroes  for
$z_L(x,-j)$ are 
easily seen to be in the algebraic closure of $\Fr\subset \bfC_\infty$ and
so have $\infty$-adic absolute value $1$.
\medskip

\begin{example} \label{trivialzeta}
In Example 3 of \cite{go2} we discussed the basic example $\zeta_A(s)$, 
$s\in S_\infty$,
of the zeta-function of $A=\Fr[\theta]$. Note obviously that $d_\infty=1$
for such an $A$ and a ``positive'' polynomial is just
a monic polynomial. Let $s=(x,y)$; one finds easily that
\begin{equation}\label{zetaform}
\zeta_A(s)=\sum_{e=0}^\infty x^{-e}\left(\sum_{\substack{n{\rm ~monic}\\
\deg(n)=e}}\langle n\rangle^{-y}\right)\,.\end{equation}
For a non-negative integer $j$ one then has
\begin{equation} \label{zetaformneg}
z_{\zeta_A}(x,-j)=\sum_{e=0}^\infty x^{-e}\left( \sum_{\substack{n{\rm
~ monic}\\\deg(n)=e}} n^j\right)\,.\end{equation}
For $e\gg 0$, the sum in parentheses vanishes (in fact, Lee shows
that one can choose $e> l_r(j)/(r-1)$ where $l_r(j)$ is the
sum of the $r$-adic digits of $j$). When $j$ is positive
and divisible by $r-1$, then $z_{\zeta_A}(x,-j)$  has a simple zero
at $x=1$. Thus
$\zeta_A(s)$ has a simple trivial zero at $s_{-j}=(\pi^j,-j)$, 
just as the Riemann zeta
function has a simple
trivial zero at negative even integers. A very similar story
happens for general $L(\rho,s)$ of Galois type.
\end{example}\smallskip

For $r$ and small positive $j$ not divisible by $r-1$, 
computer calculations have
shown that $z_{\zeta_A}(x,-j)$ is a
polynomial in $x^{-1}$ which is irreducible over
$\Fr(T)$ and has associated Galois group equal to
the full symmetric group. If $j\equiv 0\pmod{r-1}$, a similar statement is
true computationally once $1-x^{-1}$ is factored out.

For an $L$-series $L(\psi,s)$ of Drinfeld type, as well as other more general
$L$-functions, one should find trivial zeroes in
a very similar fashion. For simplicity let $\A=\Fr[T]$.
Then, for a Drinfeld module (and, more
generally, for pure, uniformizable $T$-modules), the factors at 
the infinite primes (whose zeroes are then the trivial zeroes)
should conjecturally arise in a fashion completely
analogous to the one used above for $L(\rho,s)$.
Indeed, as mentioned before, the
theory of Boeckle and Pink \cite{bp1} computes 
$$Z_L(x,-j)=\prod_{\mathfrak P~{\rm good}}P_{\mathfrak P}(n{\mathfrak P}^jx^{-\deg_{\Fr}
\mathfrak P})^{-1}$$
via the cohomology
of the crystal associated to the $\A$-module $\psi\otimes C^{\otimes j}$
(where $C$ is the Carlitz module)
through an associated trace formula. Each local 
factor $P_{\mathfrak P}(n{\mathfrak P}^j u)$ at a good prime $\mathfrak P$
may be computed via the canonical Galois action on
the Tate module $T_v(\psi \otimes C^{\otimes j})$.

Suppose that $\psi$ is defined over a finite extension $L$ of $k$.
Let $\sigma \colon L \to C_\infty$ be a
$k$-embedding and let $L_\sigma$ be the completion of $L$ under the
induced absolute value. Let $\psi^\sigma\otimes C^{\otimes j}$ 
be the $T$-module defined
over $L_\sigma$ obtained by applying $\sigma$ to the coefficients of
$\psi\otimes C^{\otimes j}$ (note that $\sigma$ acts as the identity
on the coefficients of $C$ and its tensor powers).  
Via a fundamental result of Anderson \cite{a1},
the module $\psi^\sigma\otimes C^{\otimes j}$
is {\it uniformizable} and arises from a lattice $M_{\sigma,j}$.
The action of the decomposition group at the infinite place defined
by $\sigma$ can then be computed via the Galois action on this
lattice. Therefore it must factor through a finite Galois extension 
precisely because the
lattice generates a finite extension of $L_\sigma$. Consequently the associated
characteristic polynomial of Frobenius, which should conjecturally
impart trivial zeroes to $L(\psi,s)$, will have constant (i.e., in the
algebraic closure of $\Fr\subset \bfC_\infty$) coefficients just 
as it does for $L(\rho,s)$. The product of such characteristic polynomials
over all infinite primes should then give all the trivial zeroes at $y=-j$.
A very similar description
is expected for general $\A$.

\begin{example} \label{carlitzmodule}
As mentioned above, the $L$-series $L(C,s)$ of the Carlitz module $C$ over
$\Fr(\theta)$ is $\zeta_{\Fr[\theta]}(s-1)$. From Example \ref{trivialzeta} we
see then that $L(C,s)$ has a trivial zero at $-j+1$ where $j$ runs over
the positive integers divisible by $r-1$. This is what is predicted
by the above prescription.
\end{example}

Let $v$ be a closed point in ${\rm Spec}(\A)$ and let
$L(s)$ be an $L$-series of arithmetic type. Then the logarithmic growth
of the degrees of the special polynomials also
allows one to establish the $v$-adic interpolation $L_v(x,y)$ of
the $L$-functions given above. The function $L_v(x,y)$ is naturally
defined on
the space $\bfC_v^\ast\times S_v$, where
$\bfC_v$ is the completion of an algebraic closure of the local field
$\k_v$ and $S_v$ is the completion of $\Z$ with respect to a
certain topology (one sees easily that $S_v$ is isomorphic to the 
product of a finite cyclic group $H_v$ with $\Zp$; see \S 8.3 of \cite{go1}).
We will continue to use $x$ for the first variable (now in $\bfC_v^\ast$)
and $y$ for the second variable (now in $S_v$); thus, in this case,
our notation here differs slightly from that of \cite{go2}.
These functions have analytic properties completely similar
to those possessed by the original functions on 
$S_\infty$; so we again have a $1$-parameter
family of entire power series in $x^{-1}$ with
very strong continuity properties in the variable $y$.

Notice that the zeroes of all these entire functions are algebraic
over the base completion of $\k$ (i.e., $\K$ or $\k_v$) by
standard non-Archimedean function theory.

\begin{example} \label{vadiczetafunction}
Let $\A=\Fr[T]$ and let $v$ be a prime of degree $d$. The $v$-adic
interpolation of $\zeta_A(s)$ will be denoted $\zeta_{A,v}(x,y)$.
One has
\begin{equation} \label{vadiczetasum}
\zeta_{A,v}(x,y)=\sum_{e=0}^\infty x^{-e} 
\left(\sum_{\substack{n {\rm ~monic}\\ \deg(n)=e\\ (n,v)=1 }}n^y  \right)\,,
\end{equation}
where $x\in \bfC_v^\ast$ and $y\in S_v$.
\end{example}

\begin{rem} \label{nonimpact}
Let $\rho:\Gal (L^{\rm sep}/L)\to \Aut_{\Qpbar}(V)$ be a representation
of Galois type, as above, with $L$-series $L(\rho,s)$. Let $j$ be a 
non-negative integer. It is important to note that as the trivial zeroes
of $z_L(x,-j)$ are constants, they also have $v$-adic absolute value $1$; thus
their effect $v$-adically is very limited. Under the above
conjectures on the Galois modules associated to Drinfeld modules,
etc., a similar remark should ultimately hold in complete generality.
\end{rem}

Notice that the very act of interpolating $L(\rho,s)$
$v$-adically also removes the Euler
factors at the primes above $v$ in $z_L(x,-j)$. In other words,
let 
$$z_L(v;x,-j):=
z_L(x,-j)\hat{z}_L(v;x,-j)\,$$
where $$\hat{z}_L(v;x,-j):=
\prod_{\substack{{\mathfrak P} \mid v\\{\mathfrak P}~ {\rm good}}}P_{\mathfrak P}(n{\mathfrak P}^{-(x\pi_\ast^j,-j)})\,.$$
Then $z_L(v;x,-j)=L_v(x,-j)$ where $-j\in S_v$ (and is obviously also
a polynomial in $x^{-1}$). 

\begin{defn}\label{vadictrivial}
The zeroes of $\hat{z}_L(v;x,-j)$ are the {\it $v$-adic trivial zeroes}
of $L(\rho,s)$ at $-j$.
\end{defn}

The impact of Remark \ref{nonimpact} is precisely that 
we can ignore $v$-adically 
the $\infty$-adic trivial zeroes (i.e., the trivial zeroes
of $z_L(x,-j)$) and work as above.
As usual, the union over all $j$ of these zeroes is the set of
all $v$-adic trivial zeroes. They lie in a finite 
extension of $\k_v$.

For an $L$-series $L(\psi,s)$ of Drinfeld type, the
$v$-adic trivial zeroes are given in exactly the same way. The main
difference is that the existence of the $\infty$-adic trivial-zeroes and their
$v$-adic influence is conjectural for such $L$-series at this moment.

\begin{example}\label{vadiczeta}
We continue examining the basic case of Example \ref{trivialzeta}.
Let $v$ be a finite prime of degree $d$ in $\Fr[T]$ associated to 
a monic irreducible $f(T)$ and let $j$ be a non-negative integer.
Then the $v$-adic trivial zeroes of $\zeta_A(s)$ are the elements 
$x\in \bfC_v$ with $0=1-f^jx^{-d}$; these are considered with the obvious
multiplicities when $d$ is divisible by $p$.
\end{example}

Note the remarkable similarity between
the $\infty$-adic and $v$-adic trivial zeroes. For instance, let $\zeta_A(s)$
be as in Examples \ref{trivialzeta} and \ref{vadiczeta}, and let
$v=(f)$ be a prime of degree $1$. Then the $\infty$-adic trivial
zeroes occur at $(\pi^j,-j)\in S_\infty$ for $j>0$ and divisible by $(r-1)$, while
the $v$-adic trivial zeroes occur at $(f^j,-j)\in \bfC_v^\ast \times S_v$ for non-negative $j$. Obviously both $\pi$ and $f$ are parameters in their
respective local fields.

We finish this section by using the above ideas to factor
the special polynomials. We begin at $\infty$ and let $L(s)$,
$s\in S_\infty$, be an $L$-series of arithmetic type. Let $j$ be a non-negative
integer. Then, as we have seen, $L(x,-j)$ is a polynomial in $x^{-1}$
and, conjecturally (in the case $L=L(\psi,s)$), 
there is a polynomial factorization
\begin{equation}\label{inftyfactorization}
L(x,-j)=L_{\rm triv}(x,-j)L_{\rm nontriv}(x,-j)\,,\end{equation}
where $L_{\rm triv}(x,-j)$ is the product of
the factors arising from the infinite primes and 
whose zeroes are the trivial zeroes at $-j$.
It is important to note that these polynomials may be trivial
(i.e., the constant polynomial $1$). The zeroes of $L_{\rm nontriv}(x,-j)$
are referred to as the ``non-trivial zeroes at $-j$.''
Now let $v$ be a finite prime and view $-j$ as lying in $S_v$. Let
$$L_{v,\rm triv}(x,-j)=\prod_{\substack{{\mathfrak P}\mid v\\
{\mathfrak P}~ {\rm good}}}P_{\mathfrak P}(n{\mathfrak P}^jx^{-\deg n{\mathfrak P}})\,,$$
and ``rename'' $z_L(x,-j)$ as $L_{v,{\rm nontriv}}(x,-j)$. Then by
the $v$-adic construction
we have the factorization
\begin{equation}\label{vadicfactorization}
L_v(x,-j)=L_{v,{\rm triv}}(x,-j)L_{v, \rm nontriv}(x,-j)\,.\end{equation}
The zeroes of $L_{v, \rm nontriv}(x-j)$ are then called the ``$v$-adic
non-trivial zeroes at $-j$,'' etc. Again, it is possible that these
polynomials will be identically $1$. Both the factorization
at $\infty$ and at finite primes can be put in exactly
the same form by setting $L_\infty(x,-j):=L(x,-j)$,
$L_{\infty, \rm triv}(x,-j):=L_{\rm triv}(x,-j)$, etc.

In Section \ref{counter} we will see that there is,
conjecturally, a further decomposition of these polynomials.

\begin{rem} \label{whybadfactors}
The reader may well wonder why, besides the obvious classical
analogies, one would want to have Euler factors at the finitely
many bad primes in the definition of an arithmetic $L$-series $L(s)$.
However, we have seen how removing Euler factors adds zeroes to
an $L$-series.
These zeroes might then unnecessarily enlarge the splitting field
associated to $L(s)$ and $y$ (i.e., the algebraic
extension $\K_L^{\rm tot}(y)$
of $\K_L(y)$ obtained
by adjoining the zeroes of $L(x,y)$). So, from the viewpoint of splitting
fields at least, having
such local factors is quite desirable.
\end{rem}

\section{Krasner's Lemma} \label{krasner}
In this section we recall Krasner's Lemma and put it in a form
which is particularly useful in characteristic $p$.

Let $\mathcal K$ be an arbitrary field which is
complete under a general (not necessarily discrete) non-trivial non-Archimedean
absolute value $\vert?\vert$. The
characteristic of $\mathcal K$ may also be completely arbitrary. 
Let $\overline{\mathcal K}$
be a fixed algebraic closure of $\mathcal K$ 
equipped with the canonical extension of $\vert?\vert$. Let ${\mathcal F}$ be a subfield of $\overline{\mathcal K}$ with 
 maximal separable (over
$\mathcal K$) subfield  ${\mathcal F}_s$. In particular,
$\overline{\mathcal K}_s=\mathcal{K}^{\rm sep}$= the separable closure
of $\mathcal K$ in $\overline{\mathcal K}$.

Let $\alpha\in \overline{\mathcal K}$. 

\begin{defn} \label{distance}
If $\alpha$ is not totally inseparable over $\mathcal K$ then we
set
$$\delta(\alpha)=\delta_{\mathcal K}(\alpha):=\min_{\sigma\neq id}\{\vert
\sigma (\alpha)-\alpha\vert\}\,,$$
where $\sigma$ runs over the non-identity 
$\mathcal K$-injections  of $\mathcal{K}(\alpha)$
into $\overline{\mathcal K}$. If $\alpha$ is
purely inseparable over $\mathcal K$, then we set $\delta(\alpha)=0$.
\end{defn}
\smallskip

Notice that if ${\rm char}({\mathcal K})=p>0$ then
$$\delta(\alpha^{p^i})=\delta(\alpha)^{p^i}$$ 
for $i\geq 0$.

Now let $\beta$ be another element in $\overline{\mathcal K}$.
Krasner's Lemma is then stated as follows.

\begin{prop} \label{Krasnerlemma}
Suppose that $\alpha$ is separable over ${\mathcal K}(\beta)$
and that $\vert \beta-\alpha\vert< \delta (\alpha)$. Then
${\mathcal K}(\alpha)\subseteq {\mathcal K}(\beta)$.
\end{prop}

\begin{proof}
By the separability assumption, the result follows if one knows that there are
no non-trivial embeddings of
$\mathcal{K}(\alpha,\beta)$ over $\mathcal{K}(\beta)$. But if $\tau$ is
any such injection then one has
$$\vert \tau(\alpha)-\alpha\vert =\vert (\tau(\alpha)-\beta)+(\beta -\alpha)\vert \leq \vert \beta-\alpha\vert<\delta(\alpha)$$
as $\vert \tau(\alpha)-\beta\vert=\vert \tau(\alpha-\beta)\vert=\vert \beta-\alpha\vert$. Thus $\tau=id$.\end{proof}

\begin{cor}\label{pthpowers}
Let $\alpha$ be any element in $\overline{\mathcal K}$ and suppose
that $\vert \beta-\alpha\vert<\delta(\alpha)$. Then
${\mathcal K}(\alpha)_s\subseteq {\mathcal K}(\beta)_s
\subseteq {\mathcal K} (\beta)$.
\end{cor}

\begin{proof} Suppose that $\mathcal K$ has characteristic $p>0$.
Now, for some $i\geq 0$ one knows that $\alpha^{p^i}$ is separable
over $\mathcal K$. As 
$$\vert \beta^{p^i}-\alpha^{p^i}\vert=\vert\beta -\alpha\vert^{p^i}<
\delta(\alpha)^{p^i}=\delta(\alpha^{p^i})\,,$$
the result follows from the proposition.\end{proof}

\begin{cor} \label{betaandalpha}
Suppose that $\vert \beta-\alpha\vert <\delta (\alpha)$. Then
$\delta(\beta)\leq \delta (\alpha)$ with equality if and only if
${\mathcal K}(\beta)_s={\mathcal K}(\alpha)_s$.
\end{cor}

\begin{proof} Let $\sigma$ be an injection of ${\mathcal K}(\beta)$ into
$\overline{\mathcal K}$ over $\mathcal K$. Then
$$\beta-\sigma (\beta)=( \beta - \alpha) +(\alpha -\sigma (\alpha))+
(\sigma (\alpha)-\sigma (\beta))\,.$$
The first and third terms on the right
have the same absolute value. Moreover, by
assumption, if the second term is non-zero then its absolute value
is the greatest of the three; thus it is also the absolute value of 
$\beta-\sigma(\beta)$. The result now follows. \end{proof}

Let $\mathcal K$ have characteristic $3$ and let $\lambda\in \mathcal K$
with $\vert \lambda\vert >1$. Using $\alpha:=\lambda^{1/2}$ and 
$\beta:=\alpha+\lambda^{1/{3^i}}$, for some $i>0$, 
one sees that the above corollary
cannot be strengthened to an equality between ${\mathcal K}(\alpha)$
and ${\mathcal K}(\beta)$ in general. 

Finally, the reader may trivially establish an Archimedean analogue of
Krasner's Lemma upon defining $\delta (\alpha):=\vert \alpha-\bar{\alpha}
\vert/2$ for a complex number $\alpha$.
\section{Review of some conjectures from \cite{go2}} \label{review}
Since they are used so often in this paper, we recall Conjectures
4 and 5. Let $L(s)$, $s=(x,y)\in S_\infty$, be an 
$L$-function of arithmetic type. 
We write
$$L(x,y)=\sum_{e=0}^\infty a_e(y)x^{-e}\,.$$
For each $y\in \Zp$, this power series has coefficients in the
finite extension $\K_L(y)$ of $\K$. As in
Remark \ref{whybadfactors}, we let $\K_L^{\rm tot}(y)$ be the
extension of $\K_L(y)$ obtained by adjoining the zeroes of $L(x,y)$;
we let $\K_{L,s}^{\rm tot}(y)$ be its maximal separable (over $\K_L(y)$)
subfield.

The essential part of
an algebraic extension of function fields in $1$-variable over a finite
field, whether local or
global, is the maximal separable subfield. Indeed, well-known arguments
show that totally-inseparable extensions are defined uniquely by their
degree (see Corollary 8.2.13 of \cite{go1}).
\smallskip

\noindent
{\bf Conjecture 4 of \cite{go2}}. The field
$\K_{L,s}^{\rm tot}(y)$ is a finite extension of $\K$. 
\smallskip

\noindent
The obvious $v$-adic analogue of the above conjecture is also postulated in
\cite{go2}

Viewed as power series in $x^{-1}$ for fixed $y$,
$L(x,y)$ has an associated Newton polygon in ${\mathbb R}^2$.
(To distinguish between the characteristic $p$ theory, we use
$X$ and $Y$ for the coordinates of ${\mathbb R}^2$.)

In $\bfC_\infty$ we may write
$$L(x,y)=\prod_i (1-\beta_i^{(y)}/x)\,.$$
Obviously only non-zero $\beta_i^{(y)}$ are of interest, in which case
we set $\lambda_i^{(y)}:=1/\beta_i^{(y)}$. The valuation (using $\nu_\infty$)
of 
$\lambda_i^{(y)}$, and so $\beta_i^{(y)}$, is 
computed by the Newton polygon of $L(x,y)$. Standard theory
shows that the $\beta_i^{(y)}$ tend to $0$ as $i$ tends to $\infty$,
whereas the $\lambda_i^{(y)}$ also tend to $\infty$; in fact, with a little
thought one sees that this can be made uniform with respect to $y$. We call
the $\beta_i^{(y)}$ (resp.\ $\lambda_i^{(y)}$) the ``zeroes in $x$'' of $L(s)$
(resp.\ ``zeroes in $x^{-1}$'').
An advantage of using
$\beta_i^{(y)}$ as opposed to $\lambda_i^{(y)}$ is that
the slope of a side of the Newton polygon equals the valuation of
the corresponding element $\beta_i^{(y)}$; for $\lambda_i^{(y)}$ one
needs to multiply by $-1$.
\smallskip

\noindent
{\bf Conjecture 5 of \cite{go2}}. There exists a positive real number
$b=b(y)$ such that if $\delta\geq b$, then there exists at most
one zero in $x^{-1}$ of $L(x,y)$ of absolute value $\delta$.
\smallskip

\noindent
In other words, outside of finitely many anomalous cases, zeroes are uniquely
determined by their absolute values. The conjecture is also formulated
$v$-adically.

Conjecture 5 is based on the examples of Wan, Sheats, etc.,
and appears to play a role similar to the classical Generalized Riemann
Hypothesis. Indeed in \cite{go2} we showed how it leads to a variant
of the classical Generalized Riemann Hypothesis for number fields.
It implies Conjecture 4 simply because one can
then easily show that almost all zeroes of $L(s)$ are totally inseparable
over $\K_L(y)$. To show that $\K^{\rm tot}_L(y)$ is itself a finite
extension of $\K_L(y)$ (and so of $\K$), one factors $L(x,y)$
into the $L$-series of ``simple motives'' and then applies the
Generalized Simplicity Conjecture (Conjecture 7 of \cite{go2}).

Our next section explains how to use the trivial zeroes to find 
counter-examples to Conjecture 5. We also suggest a reasonable
modification of Conjecture 5.

\section{The counter-examples}\label{counter}
Let $L(s)$, $s=(x,y)\in S_\infty$, be an arithmetic $L$-series
which we continue to write as
$\sum_{e=0}^\infty a_e(y)x^{-e}$.
Let $n$ be a positive integer and let $y_0\in \Zp$. 
Suppose that the
first $n$ slopes of the Newton polygon of $L(x,y_0)$ (as a function
of $x^{-1}$) are finite. 

\begin{lemma} \label{newtonclose}
There is an non-trivial open neighborhood $U(y_0,n)$ of
$y_0$ such that if $y\in U(y_0,n)$, then
the first $n$ segments of the Newton polygon in $x^{-1}$
of $L(x,y)$ are the same as those for $L(x,y_0)$.
\end{lemma}

\begin{proof}
The functions $a_e(y)$ are continuous. Moreover, let $\nu_\infty$
be the additive valuation associated to $\infty$. Then, from
\cite{b1}, one also has exponential
lower bounds on $\nu_\infty(a_e(y))$ which are independent of $y$.
The result follows directly.\end{proof}

\noindent
A completely analogous $v$-adic result follows in the same way.

\begin{rem} \label{equivrelations}
a. We can use  the Newton polygons  of $L(x,y)$ to define equivalence relations
on $\Zp$ (or its $v$-adic analogue $S_v=\Zp\times H_v$ where $H_v$ is a finite
abelian group etc.) in the following fashion.
Let $n$ be a fixed positive integer. Let $y_i\in \Zp$, $i=1,2$, be such that
the Newton polygon of $L(x,y_i)$ has $n$ finite slopes for each $i$.
We then say that $y_1 \sim_n y_2$ if and only if the Newton polygons
of both $L(x,y_1)$
and $L(x,y_2)$ 
have the same first $n$ segments. If $y\in \Zp$ does not have
$n$ finite slopes, then, by definition, $y$ will only be equivalent to itself.
It is clear
that $\sim_n$ is an equivalence relation which only
depends on $L(s)$ and $n$.\\
b. The impact of Lemma \ref{newtonclose} is precisely that an equivalence
class of $\sim_n$ consisting of more than one element is then open
in $\Zp$.\\
c. Let $y\in \Zp$ belong to an open equivalence class $E_y$ 
under $\sim_n$ and let
$m$ be the least non-negative integer such that $U:=y+p^m\Zp\subseteq E_y$.
Thus, on $U$, the first $n$-segments of the Newton polygon are an
invariant of the maps $z\mapsto z+\beta$ where $\beta\in p^m\Zp$. 
We believe that such statements may be viewed as possible
``micro-functional-equations'' for (the Newton polygon of) $L(x,y)$. 
See Section \ref{basic}
for an example worked out in detail.
\end{rem}

It seems reasonable that the family of Newton polygons associated
to $L$-series actually determine the $L$-series itself. We state this
more succinctly in the following question.

\begin{question}\label{langlands} Let $A=\Fr[T]$ and let $\phi_1$
and $\phi_2$ be two non-isogenous
Drinfeld modules over $\Fr(\theta)$ of the same
degree. Does the family of Newton polygons serve to distinguish
between $L(\phi_i,s)$ ($s\in S_\infty$) for $i=1,2$?
\end{question}

\noindent
Obviously, there are many variants of Question \ref{langlands} that may
also be formulated.

We can now construct the counter-examples.

\begin{example}\label{counter1}
Let $\A$ be arbitrary but where $d_\infty>1$. If $j$ is a positive
integer
divisible by $r^{d_\infty}-1$ then $\zeta_A(s)$ has  trivial
zeroes at $(\zeta\pi_\ast^j,-j)$, where $\zeta$ runs over the $d_\infty$-th
roots of $1$ with multiplicity. Thus there is a segment of the Newton
polygon of $\zeta_A(x,-j)$ (in $x^{-1}$) which has slope $j/d_\infty$ and whose
projection to the $X$-axis has length $\geq d_\infty$. Lemma
\ref{newtonclose} now assures us that all $y$ sufficiently close
to $-j$ will possess this property. We now construct a counterexample to
Conjecture 5 inductively. Let $y_0=r^{d_\infty}-1$. Let
$y_1=y_0+p^{t_1}(r^{d_\infty}-1)$ where $t_1$ is a non-negative
integer chosen (in accordance with Lemma \ref{newtonclose}) so that
the first $n$ segments of the Newton polygons at $-y_0$ and $-y_1$ are the
same and where these segments include the one associated to the
trivial zeroes at $-y_0$. Now construct $y_2$ in the same fashion
but where we choose $t_2$ to also be greater than $t_1$ etc. The
sequence $\{y_i\}$ clearly converges to a $p$-adic integer $\hat y$.
The Newton polygon in $x^{-1}$ of $\zeta_A(x,-\hat{y})$ will have
infinitely many segments whose projection to the $X$-axis will have
lengths $\geq d_\infty$. There are then two cases to discuss:\\
1. $d_\infty$ is a pure power of $p$. In this case we cannot
directly conclude that there are infinitely many zeroes of
$\zeta_A(x,-\hat{y})$ which are not uniquely determined by
their absolute value simply because we do not know a-priori
that the zeroes of $\zeta_A(x,-\hat{y})$
are not totally inseparable. However, if one
also assumes the Generalized Simplicity Conjecture 
(Conjecture 7 of \cite{go2}),
then almost all such zeroes cannot be totally inseparable 
and so Conjecture 5 must now be false.\\
2. $d_\infty$ is not a pure $p$-th power. In this case, there
are at least two distinct $d_\infty$-roots of unity. One can then choose
the $t_i$ sufficiently large so that the distinct trivial zeroes
separate the nearby zeroes. In this case, one obtains
a counter-example unconditionally.\end{example}
\medskip

One can often use Krasner's Lemma to obtain 
similar constructions as in the following example.
\medskip

\begin{example}\label{counter2}
Let $\A=\Fr[T]$. Let $f$ be a prime of degree
$d>1$ with associated place $v$ and assume that $d$
is not a pure $p$-th power. Then the $v$-adic trivial zeroes
of $\zeta_A(s)$ at $-j$ are the roots of $1-f^jx^{-d}$.
Let $j\not \equiv 0\pmod{d}$ and let $\alpha$ be one
such root. Then $\alpha\not\in \k_v$. Moreover, it is easy to
see that, in the notation of Section \ref{krasner}, we have
$$\delta_{\k_v}(\alpha)=\vert \alpha\vert_v\, .$$
Let $y\in S_v$ be sufficiently close to $-j$ so that
$\zeta_{A,v}(x,y)$ has a zero $\beta$ with
$\vert \beta-\alpha\vert_v<\vert \alpha\vert_v$.
By Corollary \ref{pthpowers}, the separable degree of
$\k_v(\beta)$ is greater than $1$. As such this
$\beta$ possesses a non-trivial Galois conjugate
$\beta^\prime$ which is also a zero of $\zeta_{A,v}(x,y)$
of the same absolute value. One can now proceed as
in Example \ref{counter1} to obtain a counter-example to
the $v$-adic version of Conjecture 5.
\end{example} \medskip

In the above example, it is easy to see that all
$v$-adic trivial zeroes belong to a finite extension of
$\k_v$. Thus, Krasner's Lemma does not allow us to deduce
a counter-example to Conjecture 4.

Conjecture 5 may still remain valid in its original
form in the much
more limited case where there exists only
one (including multiplicity!) trivial zero of a given absolute
value. Indeed, the techniques used in the above counter-examples do not
work in his case.

Ideally, one would like to negate the effects
of the trivial zeroes which permit the above counter-examples. Classically
one removes the effects of the trivial zeroes
through the use of the Gamma-factors (as in the introduction), 
which are the Euler factors
arising from the infinite primes, and the functional equation. Indeed,
the functional equation assures us that the trivial zeroes are quite far
from the critical zeroes (= all ``non-trivial'' zeroes).

In the characteristic $p$ case that we are studying, it has been known
for a long time that the Gamma-functions do not seem to be related to
the trivial zeroes of $L$-series. A philosophical explanation for this
phenomenon comes from the ``two $T$'s'' approach; indeed, $L$-series have
values in the field $\bfC_\infty$ whereas Gamma-functions, as they are 
related to exponential functions and their periods, must take values
in $C_\infty$. 

In any case, we need to find other methods in the characteristic $p$
theory.  We now sketch an approach to removing the
``harmful'' effects of trivial zeroes based on
Hensel's Lemma. The idea is simply to isolate those
zeroes which are influenced by the trivial zeroes so that they can
be removed from the conjectures and handled separately. Whether the definition
of the $L$-series should be altered, as in the classical case, to physically
remove these zeroes is unknown.

In order to isolate those zeroes which are sufficiently close to trivial
zeroes, an affirmative
answer to the following question about trivial zeroes would give the
nicest situation. This question
seems reasonable in view of examples and ramification considerations.
So let $w$ be a place of $\k$ (either $\infty$ or a finite
place) and consider the $w$-adic interpolation of an $L$-series
$L(s)$ of arithmetic type. Recall that, from Equations \ref{inftyfactorization}
and \ref{vadicfactorization}, we have a factorization
$$L_w(x,-j)=L_{w,\rm triv}(x-j)L_{w,\rm nontriv}(x,-j)\,.$$
Let $e>0$ be a real number.
Then, standard non-Archimedean analysis leads to a rational factorization
\begin{equation}\label{absolutevaluefactorization}
L_w(e;x,-j)=L_{w,\rm triv}(e;x,-j)L_{w,\rm nontriv}(e;x,-j)
\end{equation}
where $L_w(e;x,-j)$ is the product of $1-\beta/x$ where
$\beta$ runs through all zeroes in $x$
of $L_w(x,-j)$
with $\nu_w(\beta)=e$; etc. (Recall that the zeroes in $x$ of $L_w(x,y)$ 
uniformly tend to $0$ so that their valuations tend to $\infty$.)

\begin{question} \label{trivialabsolutevalue}
Let $j$ be a non-negative integer. Does there exists a constant
$C>0$ (depending only on $L(s)$ and $w$) so that for $e>C$
the polynomials $L_{w,\rm triv}(e;x-j)$ and $L_{w, \rm nontriv}(e;x,-j)$
are relatively prime polynomials in $x^{-1}$?
\end{question}

Let us assume that the above question may be answered in the
affirmative and let $e$ be as in its statement.
As $L_{w, \rm triv}(e;x,-j)$ and
$L_{w,\rm nontriv}(e;x,-j)$ are relatively prime, Hensel's Lemma
now applies to polynomials  which are close to
$L_w(e; x,-j)$ (see, eg., Thm.\ 4.1 of \cite{dgs1}). 

Now let $y$ be chosen
sufficiently close to (but {\it not} equal to)
$-j$ so that the first $m$ segments of the
Newton polygons
are the same, where $m$ is large enough so that the segment associated to 
$e$ is among the first $m$ chosen. It is reasonable
to assume that $L_w(e;x,y)$, with the obvious definition,
is also then close enough to
$L_w(e;x,-j)$ for Hensel's Lemma to apply. Thus, under these
assumptions, $L_w(e;x,y)$ will inherit a rational factorization
\begin{equation}\label{yadictrivcrit}
L_w(e;x,y)=L_{w,\rm triv}(e;x,y)L_{w, \rm nontriv}(e;x,y)\,.
\end{equation}
In other words, outside of finitely many exceptional $e$,
we would then be able to isolate those zeroes of $L_w(x,y)$ which
are influenced by the trivial zeroes. 

\begin{rem} \label{henselstill}
Even if Question \ref{trivialabsolutevalue} is answered in the negative,
one can still proceed as follows. Let $d_w(e;x,-j)=1+\cdots$ be the greatest
common divisor of $L_{w,\rm triv}(e;x,-j)$ and $L_{w,\rm nontriv}(e;x,-j)$ as
polynomials in $x^{-1}$.
We can then apply Hensel's Lemma to the relatively prime factorization
\begin{equation}\label{henselalso}
L_w(e;x,-j)=L_{w,\rm triv}(e;x,-j)d_w(e;x,-j)\times
L_{w,\rm nontriv}(e;x,-j)/d_w(e;x,-j)\,.\end{equation}
By construction, the trivial zeroes are only associated to the factor
on the left in Equation \ref{henselalso}.
\end{rem}

The zeroes of $L_{w,\rm triv}(e;x,y)$ are called the ``near-trivial
zeroes associated to $e$ and $y$,'' etc. They are precisely the zeroes
which are influenced by the original trivial zeroes. (N.B.: If $y$ is
actually a negative integer itself, there is nothing a-priori to rule
out having a near-trivial zero also
being an actual trivial zero for $y$.)
The rest of the zeroes are called the ``critical zeroes'' (in analogy
with classical theory) and these are
the ones Conjecture 5 may indeed apply to.

It remains to deal with Conjecture 4. Assuming that
Conjecture 5 is established somehow for critical zeroes, the only issue 
that remains 
is to somehow establish that the field generated by all the near-trivial zeroes
for a given $y$ is also finite over $\k_w$. 
However, the degree of $L_{w,\rm triv}(e;x,y)$ is
bounded (the example of $L(\rho,s)$ will suffice to convince the reader that
this is so). Thus it would suffice to bound the
discriminants of the maximal separable subfield of
the splitting field of $L_{w,\rm triv}(e;x,y)$
(as there are only finitely many separable extensions of a local function field
of bounded degree and discriminant, see Prop. 8.23.2 of \cite{go1}).
Needless to say, such a
problem never comes up in classical theory. However, the following
examples give some evidence in favor of such bounds in the characteristic
$p$ theory.

\begin{example} \label{f5}
Let $\A=\F_3[T]$ and let $v$ correspond to a monic prime $f$ of degree $2$.
Let $z_{\zeta_A}(x,-j)$ be as in Equation \ref{zetaformneg};
one computes easily that $z_{\zeta_A}(x,-5)=1+(T-T^3)x^{-1}$. Thus
$$z_{\zeta_A}(v;x,-5)=\zeta_{A,v}(x,-5)=(1-f^5x^{-2})(1+(T-T^3)x^{-1})\,.$$
The first factor gives the trivial zeroes and the second gives
the non-trivial zeroes.
Clearly these two factors are relatively prime. Thus for $y\in S_v$ 
sufficiently close to $-5$, Hensel's Lemma may be used. Note also that 
$f^{5/2}$ is obviously a separably algebraic element. Thus,
if 
$y$ is also close enough to $-5$ so that Krasner's Lemma applies, then
we find that the near-trivial zeroes at $y$ 
associated to $5/2$ generate $\k_v(\sqrt{f})$. Indeed,
if $\beta$ is a near-trivial zero associated to $f^{5/2}$ then
$\beta$ will also satisfy a quadratic equation over $\k_v$ (and so
we deduce equality of fields as opposed to merely
inclusion as in Lemma \ref{Krasnerlemma}). Thus
Corollary \ref{betaandalpha} implies that $\delta_{\k_v}(\beta)=
\delta_{\k_v}(f^{5/2})=\vert f^{5/2}\vert_v$.
\end{example}

\begin{example}\label{f4}
We continue with the set-up of Example \ref{f5}. One has
$z_{\zeta_A}(x,-4)=(1-x^{-1})$. Thus
$$z_{\zeta_A}(v;x,-4)=\zeta_{A,v}(x,-4)=(1-f^4x^{-2})(1-x^{-1})\,.$$
In this case, Hensel's Lemma implies that near-trivial zeroes
associated to $\pm f^2$ are in $\k_v$.
\end{example}

Simple considerations of Newton polygons imply that both $(T^3-T,-5)$
and $(1,-4)$ are critical zeroes for $\zeta_{A,v}(x,y)$. 

\section{The analytic behavior of $\zeta_{\Fp[\theta]}(s)$, $s\in S_\infty$}\label{basic}
We will use the techniques and results of Diaz-Vargas \cite{dv1}
(see also \S 8.24 of \cite{go1}) and Sheats \cite{sh1} to describe
the influence of the trivial zeroes for $\zeta_{A}(s)$, $A=\Fr[\theta]$ and
$s\in S_\infty$. We will see that, contrary to what we first
expected, all zeroes of $\zeta_{\Fp[\theta]}(s)$
are near-trivial. In fact, examples lead us to expect this to hold
for all $r$; the proof will take a detailed analysis of Sheats'
method which we hope to return to in later works.

Our first result along these lines concerns the valuation at
$\infty$ of the zeroes of $\zeta_A(s)$. Let $s=(x,y)\in S_\infty$ and write
\begin{equation}\label{zetasumform}
\zeta_A(s)=\sum_{i=0}^\infty a_i(y)x^{-i}\,. \end{equation}
As before, let $\nu_\infty$ be the normalized valuation at $\infty$
with $\nu_\infty(1/T)=1$.

\begin{prop} \label{missed}
We have $\nu_\infty(a_i(y))\equiv 0\pmod{r-1}$ for all $i$ and $y$.
\end{prop}

\begin{proof} Let $j$ be a non-negative integer and
(in the notation of \cite{sh1}) 
$$S_k^\prime(j)=\sum_{\substack{n\in \Fr[T]\\n~{\rm monic}\\\deg(n)=k}}
n^j\,.$$
The main result in \cite{sh1} is to establish a formula
originated by Carlitz for $\deg(S_k^\prime(j))$ (this
formula is then used to compute $\nu_\infty(a_i(y))$ and the
Newton polygon of $\zeta_A(x,y)$). The formula expresses
$\deg(S_k^\prime(j))$ in terms of a certain
$k+1$-tuple, called the ``greedy element,'' $(x_1,\ldots,x_{k+1})$
of non-negative integers such that $\sum x_t=j$ in such a way that there
is no carry-over of $p$-adic digits and such that the first $k$ elements
are both positive and divisible by $r-1$. From this formula one obtains
a formula for $\nu_\infty(a_i(y))$ by choosing $j$ sufficiently close to
$-y$ (see Equation 2.2 and Lemma 2.1 of \cite{sh1}). The result follows
simply by noting that this formula is linear and involves only the first
$k$-terms of the given $k+1$-tuple.
\end{proof}

\begin{cor} \label{evenval}
Let $\alpha\in \K$ be a zero of $\zeta_A(x,y)$ for $y\in \Zp$. Then
$\nu_\infty(\alpha)$ is positive and divisible by $r-1$.
\end{cor}

\begin{proof} The fact that $\nu_\infty(\alpha)$ is positive follows 
easily from general theory. To see that it is divisible by $r-1$, note that
Sheats' work shows that the Newton polygon of $\zeta_A(x,y)$
only has segments of vertical length $1$ (i.e., their projection to the
$X$-axis has unit length). Thus the divisibility follows immediately
from the proposition.
\end{proof}

We now set $r=p$ in order to use Diaz-Vargas' 
simple techniques to compute the greedy element and to avoid problems 
involving carry-over of 
$p$-adic digits (this carry-over is what makes the
general $\Fr$-case so subtle). Let $n$ be a positive integer and
let $\sim_n$ be as in Part a of Remark \ref{equivrelations}.
Our first goal is to describe explicitly the equivalence classes
of $\sim_n$ in $\Zp$. Let $y\in \Zp$ be a non-negative integer which we write
$p$-adically as $\sum_{t=0}^w c_tp^t$ where $0\leq c_t <p$ for all $t$.
We set $l(y)=l_p(y):=\sum_t c_t$ as usual.  If $y\in \Zp$ is not a non-negative
integer then we set $l(y)=\infty$.

\begin{prop} \label{diazcor}
a. Let $j$ be a non-negative integer. Then
the degree in $x^{-1}$ of $\zeta_A(x,-j)$ is $[l(j)/(p-1)]$
(where $[?]$ is the standard greatest integer function).\\
b. Let $y\in \Zp$. Then $\zeta_A(x,y)$ has at least $n$ distinct
slopes if and only if $[l(-y)/(p-1)]\geq n$.
\end{prop}

\begin{proof} The first part follows immediately from Diaz-Vargas'
construction of the greedy element (see e.g., the proof of
Lemma 8.24.11 of \cite{go1}). The second part follows from the
first part and the fact that all segments of the Newton polygon
of $\zeta_A(x,y)$ are known to have projections to the $X$-axis
of unit length.\end{proof}

Note that Part a of the proposition allows one to compute explicitly
the zeroes in $y$ of $a_i(y)$ (as defined in Equation \ref{zetasumform}).

Now let $y\in \Zp$ be chosen so that $[l(-y)/(p-1)]\geq n$ and expand
$-y$ $p$-adically as $\sum_{t=0}^\infty c_i p^i$ (where it may happen
that all but
finitely many of the $c_i$ vanish). Set
$$y_n=\sum_{i=0}^e c_ip^i$$
where $\sum_{i=0}^e c_i=n(p-1)$ and $c_e\neq 0$.
Clearly $y_n\equiv 0\pmod{p-1}$.

\begin{prop} \label{diazcor2}
a. We have $-y_n\sim_n y$.\\
b. $y_n$ is the smallest element in the set of positive integers $i$
with $-i\sim_n y$.\\
c. Let $y$ and $z$ be in $\Zp$. Then $y\sim_n z$ if and only
if $y_n=z_n$.
\end{prop}

\begin{proof} This again follows from Diaz-Vargas' construction
of the greedy element.\end{proof}

Thus the open equivalence classes under $\sim_n$ are in one to one
correspondence with negative integers $-j$ with $l(j)=n(p-1)$ (and, in 
particular, $j$ is divisible by $p-1$). Let $E$ be the equivalence
class of one such $-j$ and write $j$ $p$-adically
as $\sum_{t=0}^u c_tp^t$ where $c_u\neq0$. Let $\beta\in p^{u+1}\Zp$. It is
then clear that $E$ is stable under the mapping $z\mapsto z+\beta$.

We finish this section by reworking the above results in a way
which makes more transparent the close connection all zeta
zeroes (at $\infty$!) have with trivial zeroes.
Thus let $y\in \Zp$ be arbitrary and let $\alpha$ be a zero (in
$x$) of $\zeta_A(x,y)$. From Corollary \ref{evenval} we know
that $\nu_\infty(\alpha)$ is both positive and divisible by $p-1=r-1$.
Set $j:=\nu_\infty(\alpha)$.

\begin{prop} \label{jandn}
Let $n(j,y)$ be the number of zeroes $\beta$ of $\zeta_A(x,y)$ 
with $\nu_\infty(\beta)\leq j$.\\
a. We have $n(j,y)=l(j)/(p-1)$.\\
b. We have $-j \sim_n y$.\\
c. The zero of $\zeta_A(x,-j)$ corresponding to $\alpha$ is precisely
the trivial zero of $\zeta_A(x,-j)$.
\end{prop}

\begin{proof} Clearly the trivial zero of
$\zeta_A(x,-j)$ has valuation $j$ and it is easy to see that this
is the unique zero of $\zeta_A(x,-y)$ of highest valuation. The result
now follows as before.\end{proof}

\begin{cor} \label{log}
We have $n(j,y)=O(\log (j))\,.$
\end{cor}

Let $y=-1$. The $i$-th slope of the Newton polygon
of $\zeta_A(x,-1)$ is $p^{i+1}-1$ and it is easy to see that
$n(p^{i+1}-1,-1)=i$ is asymptotic to $\log_p(p^{i+1}-1)$. Thus the number of
zeroes of $\zeta_A(x,-1)$ of valuation $\leq x$, for a positive real $x$,
is asymptotic to $\log_p(x)$. Of course many other such examples may
be worked out.

In any case, one sees that all zeroes of $\zeta_{\Fp[\theta]}(s)$
are near-trivial. For general $\Fr[T]$, calculations indicate that
Parts b and c of Proposition \ref{jandn} should
remain valid. If so, then Part c of Proposition \ref{jandn}
may ultimately afford an explanation why the results of Wan, Diaz-Vargas,
Thakur, Poonen and Sheats were obtainable by elementary means. Moreover, it
also shows that, as of this writing, we have had precious little
experience with critical zeroes.

It is also reasonable to expect that Corollary \ref{log} will be true for
all arithmetic $L$-series at all primes. A much more interesting question
is whether some version of Part a of Proposition \ref{jandn} will be
true. That is, is the analogue of $n(j,y)$ independent of $y$?

\section{Taylor expansions of classical $L$-series}\label{taylor}
Let $\A=\Fr[T]$ and consider $\zeta_A(s)$, $s=(x,y)\in S_\infty$
as in Example \ref{trivialzeta}.
It is clear from Equation \ref{zetaform} that for all $y\in \Zp$, 
$\zeta_A(x,y)$ is a power series in $x^{-1}$ with coefficients in
$\K=\k_\infty$. In this section we will establish in
great generality a very similar result
for the complex analytic
functions $\Xi(\chi,t)$ of \cite{go2} (the definition of $\Xi(\chi,t)$
will be recalled below). Consequently, as mentioned in the
introduction, these Taylor expansions
reflect the (conjectured!) rationality of their zeroes in a simpler
fashion than one finds for arbitrary entire complex functions.

\begin{defn} Let $p(t)=\sum c_jt^j$ be a non-zero complex power series.
We say that $p(t)$ is {\it almost real} if and only if
$$p(t)=\alpha h(t)$$
where $\alpha\in \mathbb{C}^\ast$ and where $h(t)$ is a non-zero
power series with real coefficients.
\end{defn}

\begin{prop} \label{galfe}
A complex power series $p(t)=\sum c_j t^j$ is almost real if and only if
the coefficients $c_j$ satisfy the ``Galois functional equation''
$$\overline{c_j}=w c_j$$
for a fixed complex number $w$ of absolute value $1$.
\end{prop}

\begin{proof}
Suppose that $p(t)=\alpha h(t)$ is almost real, where $\alpha$
is non-zero and $h(t)\in \mathbb{R}[[t]]$. Put $w:=\overline{\alpha}/
\alpha$; it is simple to check that with this $w$ the Galois functional
equation holds. Conversely, assume the Galois functional equation and
let $j_1$ and $j_2$ be two non-negative integers such that
$c_{j_1}\neq 0$. Then
$$\overline{c_{j_2}/c_{j_1}}=\overline{c_{j_2}}/\overline{c_{j_1}}=(w c_{j_2})/
(w c_{j_1})=c_{j_2}/c_{j_1}\,;$$
thus $c_{j_2}/c_{j_1}$ is real. Now let $j_0$ be the smallest
non-negative integer with $c_{j_0}\neq 0$. Then 
$$p(t)=c_{j_0}\times t^{j_0}(1+\sum_{i=1}^\infty b_it^i)$$
with $b_i$ real, and the result is established. \end{proof}

Now let $\chi$ be a non-trivial finite abelian character associated to
a Galois extension of number fields $L/k$. Let $L(\chi,s)$ be the
classical (complex) $L$-series and let $\Lambda (\chi, s)$ be the completed
$L$-function with the Euler factors at the infinite primes. 
As is standard $\Lambda (\chi,s)$ is entire and there is a functional
equation connecting $\Lambda (\chi,s)$ and $\Lambda (\overline{\chi},1-s)$.
In particular,
$$\Lambda (\overline{\chi},1-s) =w(\chi)\Lambda (\chi, s) \,,$$
where $w(\chi)$ has absolute value $1$.
We then set $\Xi(\chi,t):=\Lambda(\chi, 1/2+it)$ following Riemann. Let
$$\Xi(\chi,t)=\sum_{n=0}^\infty a_n t^n\,,$$
be the Taylor expansion of $\Xi(\chi,t)$ about the origin.

\begin{prop} \label{anfe}
We have
$$\Lambda (\overline{\chi},1-s) =w(\chi)\Lambda (\chi,s)$$
if and only if the coefficients $\{a_n\}$ satisfy 
$$\overline{a_n}=w(\chi) a_n\,,$$
for all $n$.
\end{prop}

\begin{proof}

We know that $\overline{\Lambda (\chi,s)}=
\Lambda (\overline{\chi}, \overline{s})$. Thus we see
$$\overline{\Xi(\chi,t)}=\overline{\Lambda (\chi,1/2+it)}=
\Lambda (\overline{\chi},1/2-i\overline{t})=
\Xi(\overline{\chi},-\overline{t})\,;$$
consequently, $\overline{\Xi(\chi,\overline{t})}=\Xi(\overline{\chi},-t)$.
On the other hand, the functional equation immediately gives us
\begin{eqnarray*}
\Xi(\overline{\chi},-t)&=&\Lambda(\overline{\chi},1/2-it)\\
&=&\Lambda(\overline{\chi},1-(1/2+it))\\
&=&w(\chi)\Lambda(\chi,1/2+it)=w(\chi) \Xi(\chi,t)\,.\end{eqnarray*}
Consequently we deduce that
\begin{equation} \label{funceq}
\overline{\Xi(\chi,\overline{t})}=w(\chi)\Xi(\chi,t)\,.\end{equation}
The only if part now follows upon substituting in the power series 
for $\Xi(\chi,t)$. The if part follows since these
calculations are reversible.
\end{proof}

\begin{theorem} \label{main}
The existence of a classical style functional equation
for $\Lambda (\chi,s)$ is equivalent to the Taylor expansion
at the origin $t=0$ 
of $\Lambda (\chi,1/2+it)$ being an almost real
power series.
\end{theorem}

\begin{proof} This follows directly from Propositions
\ref{anfe} and \ref{galfe}.\end{proof}

For Dedekind zeta functions a completely similar result may easily
be established
along with some vanishing of the Taylor coefficients.
In many instances it is known that classical $L$-series may be factored
as infinite products over their zeroes. Such a factorization gives
another approach to showing that the Taylor expansion of $\Xi(\chi,t)$ is
almost real.

\end{document}